\documentclass[12pt]{amsart}
\usepackage{amssymb,amsmath}

\headsep=1cm \numberwithin{equation}{section}
\newtheorem{theorem}{Theorem}[section]
\newtheorem{lemma}[theorem]{Lemma}
\newtheorem{proposition}[theorem]{Proposition}
\newtheorem{corollary}[theorem]{Corollary}

\theoremstyle{definition}

\theoremstyle{remark}
\newtheorem{remark}[theorem]{Remark}
\newtheorem{example}[theorem]{Example}

\newtheorem{acknowledgement}{Acknowledgement}

\newcommand{\grade}{\operatorname{grade}}

\newcommand{\depth}{\operatorname{depth}}

\newcommand{\reg}{\operatorname{reg}}

\newcommand{\fm}{\frak{m}}

\newcommand{\fa}{\frak{a}}

\begin{document}
\dedicatory{Dedicated to Professor Maria Evelina Rossi}
\author[Mafi and Naderi]{ Amir Mafi and Dler Naderi}
\title[On the Hilbert coefficients...]
{On the Hilbert coefficients, depth of associated graded rings and reduction numbers}

\address{A. Mafi, Department of Mathematics, University of Kurdistan, P.O. Box: 416, Sanandaj,
Iran.} \email{a\_mafi@ipm.ir}
\address{D. Naderi, Department of Mathematics, University of Kurdistan, P.O. Box: 416, Sanandaj,
Iran.}
\email{dler.naderi65@gmail.com}


\begin{abstract} 
Let $(R,\fm)$ be a $d$-dimensional Cohen-Macaulay local ring, $I$ an $\fm$-primary ideal of $R$ and $J=(x_1,...,x_d)$ a minimal reduction of $I$.
We show that if $J_{d-1}=(x_1,...,x_{d-1})$ and $\sum\limits_{n=1}^\infty\lambda{({{I^{n+1}}\cap J_{d-1}})/({J{I^n} \cap J_{d-1}})=i}$ where i=0,1, then $\depth G(I)\geq{d-i-1}$. Moreover, we prove that if
${e_2}(I)=\sum\limits_{n = 2}^\infty{(n-1)\lambda({{I^{n }}}/{J{I^{n-1}}})}-2;$ or if $I$ is integrally closed and\\ ${e_2}(I)=\sum\limits_{n=2}^\infty{(n-1)\lambda({{I^{n}}}/{J{I^{n-1}}})}-3$, then ${e_1}(I)=\sum\limits_{n=1}^\infty{\lambda({{I^{n }}}/{J{I^{n-1}}})}-1,$ where the integers $e_i$ are the Hilbert coefficients of $I$. In addition, if $J$ is a minimal reduction of $I$ then we prove that the reduction number $r_J(I)$ is independent of $J$.
\end{abstract}

\maketitle

\section{Introduction}
Throughout the paper we will assume that $(R,\fm)$ is a $d$-dimensional Cohen-Macaulay local ring having an infinite residue field and $I$ is an $\fm$-primary ideal of $R$. An ideal $J\subseteq I$ is called a reduction of $I$ if $I^{r+1}=JI^r$ for some nonnegative integer $r$ (see \cite{Nr}). The least such $r$ is called the reduction number of $I$ with respect to $J$ and is denoted by $r_J(I)$. A reduction $J$ is called a minimal reduction if it does not properly contain a reduction of $I$. Under our assumption every minimal reduction is generated by a regular sequence. The reduction number of $I$ is defined as $r(I)=\min\{r_J(I): J$ is a minimal reduction of $I\}$. The reduction number $r(I)$ is said to be independent if $r(I)=r_J(I)$ for all minimal reductions $J$ of $I$.
Sally in \cite{S1} raised the following question: If $(R,\fm)$ is a $d$-dimensional Cohen-Macaulay local ring having an infinite residue field, then is $r(\fm)$
independent? A natural extension of this question is to replace $r(\fm)$ with $r(I)$. Let $G(I)=\bigoplus_{n\geq 0}I^n/I^{n+1}$ be the associated graded ring of $I$. Huckaba in \cite{H} and Trung in \cite{T} independently proved that if $\depth G(I)\geq{d-1}$, then $r(I)$ is independent(see also \cite{M}, \cite{Ho}, \cite{Ho1} and \cite{St}). Moreover, Wu in \cite{Wu} with some conditions proved that if $\depth G(I)\geq{d-2}$, then $r(I)$ is independent.
However, if $d\geq 2$ and $\depth G(I)\leq{d-2}$, then $r(I)$ is not independent in general. Counter-examples have been obtained in \cite{H}, \cite{M} and \cite{M1}.

The Hilbert function of $I$ is  the numerical function $H_I(n)=\lambda(R/I^n)$ (where $\lambda(.)$ denotes length) that measures the growth of the length of $R/I^n$ for all $n\geq 1$.  It is well known that for $n\gg 0$, $H_I(n)$ is a
polynomial in $n$ of degree $d$  $$ P_I(n)=\sum_{i=0}^d(-1)^ie_i(I){{n+d-i-1}\choose{d-i}},$$ called the Hilbert polynomial of $I$,  whose coefficients $e_0(I), e_1(I),..., e_d(I)$ are uniquely determined by $I$ and called the Hilbert coefficients of $I$.

Valabrega and Valla in \cite{Vv} obtained that $G(I)$ is Cohen-Macaulay if and only if there exists a minimal reduction $J$ of $I$ such that $I^n\cap J=I^{n-1}J$ for all $n$. Later on, Guerrieri in \cite{G1} asked that if $J$ is a minimal reduction of $I$ such that $\sum_{n\geq 1}\lambda({I^n\cap J}/{I^{n-1}J})=t$, then
is $\depth G(I)={d-t}$? The case $t=0$ is simply a restatement of the Valabrega-Vall theorem, whereas the case $t=1$ was proved in \cite{G1}. Some partial answers in the cases $t=2$ or $t=3$ were also proved in \cite{G2} and \cite{Gr}.
 Huckaba and Marley in \cite{Hm} and Vaz Pinto in \cite{Va} independently showed that $e_1(I)\leq\sum_{n\geq 1}\lambda({I^n}/{I^{n-1}J})$ and equality holds for some minimal reduction $J$ of $I$ if and only if $\depth G(I)\geq{d-1}$.
Another closely related conjecture was raised by Wang in \cite{W} while attempting to prove Guerrieri's conjecture. Namely, he asked whether the difference $\sum_{n\geq 1}\lambda({I^n}/{I^{n-1}J})-e_1(I)=s\geq 0$, implies $\depth G(I)\geq{d-s-1}$.
 Wang first showed that an affirmative answer to his conjecture implies the validity of Guerrieri's conjecture. Then in \cite{W} he settled the case $s=1$ (see also \cite{P}). Unfortunately, both conjectures fail in general as shown in \cite{W1}.

 Corso, Polini and Rossi in \cite{Cpr} established a general upper bound on $e_2(I)$, which is reminiscent of the bound on $e_1(I)$ due to Huckaba and Marley in \cite{Hm} and Vaz Pinto in \cite{Va}.
 Namely, it holds that $e_2(I)\leq\sum_{n\geq 2}{(n-1)}\lambda({I^n}/{I^{n-1}J})$ for any minimal reduction reduction $J$ of $I$. Furthermore, the upper bound is attained if and only if $\depth G(I)\geq{d-1}$. In addition, if $e_2(I)\geq\sum_{n\geq 2}{(n-1)}\lambda({I^n}/{I^{n-1}J})-2$ or if $I$ is integrally closed and $e_2(I)\geq\sum_{n\geq 2}{(n-1)}\lambda({I^n}/{I^{n-1}J})-4$, then $\depth G(I)\geq{d-2}$.

 In this paper we prove the following results.

 \begin{theorem}\label{th1} Let $J=(x_1,...,x_d)$ be a minimal reduction of $I$.\\
 (1) Suppose that one of the following
conditions holds:
\begin{itemize}
\item[(i)]
 ${e_2}(I)=\sum\limits_{n = 2}^\infty  {(n - 1)\lambda({{I^{n }}}/{J{I^{n-1}}})}-2; $
 \item[(ii)]
 $I$ is integrally closed and ${e_2}(I)=\sum\limits_{n = 2}^\infty  {(n - 1)\lambda({{I^{n }}}/{J{I^{n-1}}})}-3$.
\end{itemize}
Then  ${e_1}(I)=\sum\limits_{n = 1}^\infty  {\lambda({{I^{n }}}/{J{I^{n-1}}})}-1.$
 Moreover, we give a counterexample showing that   ${e_1}(I)=\sum\limits_{n = 1}^\infty  {\lambda({{I^{n }}}/{J{I^{n-1}}})}-1$ but the conditions (i) and (ii) do not hold.\\
 (2) If $J_{d-1}=(x_1,...,x_{d-1})$ and $\sum\limits_{n=1}^\infty\lambda{({{I^{n+1}} \cap J_{d-1}})/({J{I^n} \cap J_{d-1}})=i}$ where $i=0$ or $i=1$, then $\depth G(I)\geq{d-i-1}$.

 \end{theorem}

 \begin{theorem}

 Let $J=(x_1,...,x_d)$ be a minimal reduction of $I$.\\
 (1) Suppose that one of the following
conditions holds:
\begin{itemize}
\item[(i)]
 ${e_2}(I)=\sum\limits_{n = 2}^\infty  {(n - 1)\lambda({{I^{n }}}/{J{I^{n-1}}})}-2; $
 \item[(ii)]
 $I$ is integrally closed and ${e_2}(I)=\sum\limits_{n = 2}^\infty  {(n - 1)\lambda({{I^{n }}}/{J{I^{n-1}}})}-3$.
 \end{itemize}
 Then $r(I)$ is independent. Also, we study the independence of $r(I)$ with some other conditions.
  \end{theorem}

 For any unexplained notation or terminology, we refer the reader to \cite{Bh} and \cite{Rv}.

\section{ Preliminaries}
In this section we recall some known results which is studied in \cite{Hm}.
An element $x\in I\setminus I^2$ is said to be superficial for $I$ if there is an integer $c$ such that $(I^{n+1}:x)\cap I^c=I^n$ for all $n\geq c$. If $\grade I\geq 1$ and $x$ is a superficial element, then $x$ is a regular element of $R$ and by Artin-Rees Theorem $I^{n+1}:x=I^n$ for all $n$ sufficiently large.
If $R/{\fm}$ is infinite, then a superficial element for $I$ always exists. A sequence $x_1,...,x_s$ is called a superficial sequence for $I$ if $x_1$ is superficial for $I$ and $x_i$ is superficial for $I/(x_1,...,x_{i-1})$ for $2\leq i\leq s$. If $J$ is a minimal reduction of $I$, then there is a superficial sequence $x_1,...,x_d$ in $I$ such that $J=(x_1,...,x_d)$. For any element $x\in I$ we let $x^{*}$ denote the image of $x$ in $I/I^2$. We note that if $x^{*}$
is a regular element of $G(I)$, then $x$ is a regular element of $R$ and $G(I/{(x)})\cong G(I)/{(x^{*})}$ (see \cite{Hm}).

A set of ideals $\mathcal{F}={\left\{ {{I_n}}\right\}_{n \in\mathbb{Z}}}$ of $R$ where $I_n=R$ for all $n\leq 0$ and $I_1=I$ is called a {\it Hilbert filtration} if we have $(i)$ $I_{n+1}\subseteq I_n$ for all $n\geq 0$, $(ii)$ $I_nI_m\subseteq I_{n+m}$ for all $n,m\geq 0$, and $(iii)$ there is a $k\geq 0$ such that $I^n\subseteq I_n\subseteq I^{n-k}$ for all $n\geq 0$.
An element $x\in I_1$ is called superficial for $\mathcal{F}$ if there exists an integer $c$ such that $(I_{n+1}:x)\cap I_c=I_n$ for all $n\geq c$. A sequence $x_1,\ldots, x_l$ is called a superficial sequence for $\mathcal{F}$ if $x_1$ is superficial for $\mathcal{F}$ and $x_i$ is superficial for $\mathcal{F}/(x_1,\ldots, x_{i-1})$ for $2\leq i\leq l$ (see \cite{Hm}).
Let $\mathcal{F}$ be a Hilbert filtration and
$\underline x={x_1},...,{x_l}\in {I_1}$.  For any integer $n$, Huckaba and Marley in \cite{Hm} constructed the complex $C_.(\underline x,\mathcal{F},n)$ as follows: for $l=1$ they define $C_.(x_1,\mathcal{F},n)$ to be the complex $$0\longrightarrow R/I_{n-1}\overset{x_1}\longrightarrow R/I_n\longrightarrow 0.$$  For $l>1$, assume that $C_.(x_1,\ldots,x_{l-1},\mathcal{F},n)$ has been defined and consider the chain map $f: C_.(x_1,\ldots,x_{l-1},\mathcal{F},n-1)\longrightarrow C_.(x_1,\ldots,x_{l-1},\mathcal{F},n)$ given by multiplication by $x_l$. Define
$C_.(x_1,\ldots,x_{l-1},\mathcal{F},n)$ to be the mapping cylinder of $f$. Thus one can show that $C_.(x_1,\ldots,x_{l-1},\mathcal{F},n)$
has the following form
$$0\longrightarrow R/{I_{n-l}}\longrightarrow (R/{I_{n-l+1}})^l\longrightarrow (R/{I_{n-l+2}})^{{l\choose 2}}
\longrightarrow \quad ...\quad \longrightarrow R/{I_n}\longrightarrow 0.$$
Let ${C_.}(n) = {C_.}({x_1},{x_2},...,{x_l},\mathcal{F},n)$ and ${C_.}^\prime (n) = {C_.}({x_1},{x_2},...,{x_{l - 1}},\mathcal{F},n)$. For any $n$ there is an
exact sequence of complexes
$$0\longrightarrow {C_.}^\prime(n)\longrightarrow {C_.}(n)\longrightarrow{C_.}^\prime (n-1)[- 1]\longrightarrow 0.$$
Thus, we have the corresponding long exact sequence on homology:
$$\cdots\longrightarrow {H_i}({C_.}^\prime(n))\longrightarrow {H_i}({C_.}(n))\longrightarrow {H_{i-1}}({C_.}^\prime (n - 1)){\kern 1 pt}\; \stackrel{x_l}{\longrightarrow} {H_{i - 1}}({C_.}^\prime (n)) \longrightarrow\cdots. (*)$$

 Since $\mathcal{F}$ is a Hilbert filtration, ${H_i}({C_.}(\underline x,\mathcal{F},n))$ has finite length for all $i$ and
$n$. For  $i\ge 1$, consider\\
$${h_i}(\underline x,\mathcal{F}):=\sum\limits_{n = 1}^\infty{\lambda({H_i}({C_.}(\underline x,\mathcal{F},n)))}$$
 and
 $${k_i}(\underline x,\mathcal{F}):=\sum\limits_{n=2}^\infty{(n-1)\lambda({H_i}({C_.}(\underline x, \mathcal{F},n)))}.$$

These integers are well-defined by \cite[Lemma 3.6]{Hm}. Also note that, ${h_i}(\underline x,\mathcal{F})={k_i}(\underline x,\mathcal{F})=0$ for all $i > l$.
For ${\underline x^\prime}= {x_1},...,{x_{l-1}}$  we define
$${h'_i}(\underline x',\mathcal{F}):=\sum\limits_{n=1}^\infty{\lambda({H_i}({C_.}(\underline x^\prime,\mathcal{F}, n)))}.$$

\section{The associated graded ring and the first and second Hilbert coefficients}
\begin{lemma}\label{3.1}
(\cite[Theorem 3.7]{Hm}) Let $\mathcal{F}$ be a Hilbert filtration and $\underline x={x_1},...,{x_l}$ a regular sequence on R and a superficial sequence for $\mathcal{F}$. Then for each $i\ge 1$
$\sum\limits_{j\ge i}^{} {{{( - 1)}^{j - i}}{h_j}(\underline x},\mathcal{F})\ge 0$.
Moreover, equality occurs if and only if $\grade(\underline x^{*})\ge{l-i+1}$.
\end{lemma}

\begin{lemma}
Let $\mathcal{F}$ be a Hilbert filtration and $\underline x={x_1},...,{x_l}$ be a regular sequence on $R$ and a superficial sequence for $\mathcal{F}$. Then, for each $i\ge 1$,
$\sum\limits_{j\ge i+1}^{} {{{( - 1)}^{j-i-1}}{k_j}(\underline x},\mathcal{F})\ge 0$.
Moreover, equality occurs if and only if  $\grade(\underline x^{*})\ge{l-i}$.
\end{lemma}

\begin{proof}
 Fix  $i\ge 1$ and for each $n$ let $B_n$ be the kernel of the map
 ${H_i}({C_.}(n)) \longrightarrow {H_{i-1}}({C_.}^\prime (n-1))$
 given in $(*)$. Then, for each $n$, we have the exact sequence
$$0\longrightarrow {H_l}({C_.}(n))\longrightarrow...\longrightarrow {H_i}({C_.}^\prime (n-1))\overset{x_l}{\longrightarrow}{H_i}({C_.}^\prime(n)) \longrightarrow {B_n}\longrightarrow 0.$$
Therefore, for each $n$, we have
$$\lambda({B_n})=\sum\limits_{j=i+1}^l{{{(-1)}^{j-i-1}}}\lambda ({H_j}({C_.}(n)))+\sum\limits_{j=i}^l{{{( - 1)}^{j - i}}\Delta[\lambda ({H_j}({C_.}^\prime
} (n)))]$$
 and also by using \cite[Lemma 2.7]{H1} we have
$$\sum\limits_{n=2}^\infty{(n-1)\Delta[\lambda ({H_j}({{C'}_.}(n)))]}=-\sum\limits_{n=1}^\infty{\lambda ({H_j}({{C'}_.}(n)))=-{h'_j(\underline x^{\prime},\mathcal{F})}}.$$
Thus we see that
$$\sum\limits_{n =2}^\infty{(n -1)\lambda ({B_n})=\sum\limits_{j =i+1}^l{{{(-1)}^{j-i-1}}{k_j(\underline x ,\mathcal{F})}}}-\sum\limits_{j=i}^l{{{(-1)}^{j - i}}{{h'}_j(\underline x^{\prime},\mathcal{F})}}.\ \ (**)$$
By Lemma 3.1, we have
$\sum\limits_{j = i}^l {{{(-1)}^{j-i}}{{h'}_j(\underline x^{\prime},\mathcal{F})}} \ge 0$
and so
$\sum\limits_{j=i+1}^l {{{( - 1)}^{j-i-1}}{k_j(\underline x,\mathcal{F})}}\ge 0$
for each $i\ge1$.
\\
By \cite[Proposition 3.3]{Hm}, if $\grade(\underline x^{*})\ge{l-i}$ then $H_j( {C_.}(n))=0$ for all $n$ and $j\ge i+1$.
Thus, $k_j(\underline x,\mathcal{F})=0$ for $j\ge i+1$.
\\
Conversely, suppose for $i\ge 1$
$$\sum\limits_{j=i+1}^l {{{(-1)}^{j-i-1}}{k_j(\underline x,\mathcal{F})}}=0.$$
So by $(**)$, $\sum\limits_{n=2}^{\infty} {(n-1)\lambda({B_n})}=0$. Then $\sum\limits_{n=1}^{\infty} {\lambda({B_n})}=\sum\limits_{j=i+1}^l {{{(-1)}^{j-i-1}}{h_j(\underline x,\mathcal{F})}}=0$ and by Lemma \ref{3.1} we obtain $\grade(\underline x^{*})\ge {l-(i+1)+1=l-i}.$
\end{proof}

Following Marley \cite{M1}, given a function $f: \mathbb{Z}\longrightarrow\mathbb{Z}$, define the first difference function of $f$, $\Delta^1(f):\mathbb{Z}\longrightarrow\mathbb{Z}$ by $\Delta^1(f)(n)=f(n+1)-f(n)$. We usually write $\Delta^1(f(n))$ for $\Delta^1(f)(n)$. Inductively define the ith difference function of $f$, $\Delta^i(f):\mathbb{Z}\longrightarrow\mathbb{Z}$ , by $\Delta^i(f(n))=\Delta^1(\Delta^{i-1}(f(n)))$.
Following Huckaba and Marley \cite{Hm}, a reduction of a filtration $\mathcal{F}$ is an ideal $J\subseteq I$ such that $JI_n=I_{n+1}$ for all large $n$. A minimal reduction of $\mathcal{F}$ is a reduction of $\mathcal{F}$ minimal with respect to containment.
\begin{remark}
By \cite[\S 4.]{Hm} we have

$${\Delta^d}[{P_\mathcal{F}}(n)-{H_\mathcal{F}}(n)]=\lambda({{I_n}}/{J{I_{n-1}}})-\sum\limits_{i=2}^d {{{(-1)}^i}\lambda({H_i}({C_.}(n)))}$$
and
$${e_i}(\mathcal{F})=\sum\limits_{n=i}^\infty{\left({{}_{i-1}^{n-1}} \right)}\;{\Delta ^d}[{P_\mathcal{F}}(n)-{H_\mathcal{F}}(n)],$$ where $J$ is a minimal reduction of $\mathcal{F}$ and ${e_i}(\mathcal{F})$ is the Hilbert coefficients of $\mathcal{F}$.
Therefore we obtain the following
$${e_1}(\mathcal{F})=\sum\limits_{n=1}^\infty{\lambda ({{I_n}}/{J{I_{n - 1}}})}-\sum\limits_{i = 2}^d {{{(-1)}^i} {h_i}(\underline x,\mathcal{F})}$$
and
$${e_2}(\mathcal{F})=\sum\limits_{n=2}^\infty{(n-1)\lambda ({{I_n}}/{J{I_{n - 1}}})}-\sum\limits_{i=2}^d {{{( - 1)}^i}{k_i}(\underline x, \mathcal{F})}.$$
\end{remark}

\begin{proposition} (compare with \cite[Theorem 3.1]{Cpr})
 Let $\mathcal{F}$ be a Hilbert filtration and $J=(x_1,...,x_d)$ be a minimal reduction of $\mathcal{F}$. Then
\[{e_2}(\mathcal{F})\le\sum\limits_{n=2}^\infty{(n-1)\lambda({{I_{n }}}/{J{I_{n-1}}})}.\] The equality holds if and only if  $\depth G(\mathcal{F})\ge{d-1}.$
\end{proposition}
\begin{proof}
By Remark 3.3 and Lemma 3.2 we have the following
$${e_2}(\mathcal{F})=\sum\limits_{n = 2}^\infty{(n - 1)\lambda({{I_{n }}}/{J{I_{n-1}}})}-\sum\limits_{i=2}^d {{{(-1)}^{i-2}}{k_i(J,\mathcal{F})}},$$ and
 $$\sum\limits_{i=2}^d {{{(-1)}^{i-2}}{k_i(J,\mathcal{F})}}\ge 0.$$ Thus
\[{e_2}(\mathcal{F})\le\sum\limits_{n=2}^\infty{(n-1)\lambda({{I_{n }}}/{J{I_{n-1}}})}.\]
 Also the equality follows by Lemma 3.2.
\end{proof}

\begin{remark}\label{3.5}
By \cite[Lemma 3.2]{Hm} and $(*)$ we obtain the following exact sequence

$$0\longrightarrow {({I_{n-l+1}}:(\underline{x}))}/{{I_{n - l}}}\longrightarrow{({I_{n - l + 1}}:(\underline{x}'))}/{{I_{n - l}}}\overset  {x_l}{\longrightarrow}{({I_{n -l +2}}:(\underline{x}'))}/{{I_{n-l+ 1}}}\longrightarrow ...$$

$$\longrightarrow {({I_{n-1}}\cap (\underline{x}'))}/{{(\underline{x}'){I_{n-2}}}}\overset{x_l}{\longrightarrow}{({I_n}\cap (\underline{x}'))}/{(\underline{x}'){I_{n -1}}} \longrightarrow {({I_{n}}\cap(\underline{x}))}/{(\underline{x}){I_{n -1}}}\longrightarrow$$

 $$ R/{({I_{n-1}},(\underline{x}'))}\overset{x_l}{\longrightarrow} R/{({I_n},(\underline{x}'))} \longrightarrow R/{({I_n},(\underline{x}))} \longrightarrow 0.$$

If $A_{n}$ is the kernel of the map
$R/{({I_{n - 1}},(\underline{x}'))}\overset{x_l}\longrightarrow R/{({I_n},(\underline{x}'))}$,
then
$A_{n}=(({I_n},(\underline{x}')):{x_l})/{({I_{n-1}},(\underline{x}'))}$.\\
If $C_{n}$ is the kernel of the map
${({I_{n }} \cap (\underline{x}))}/{(\underline{x}){I_{n - 1}}}\longrightarrow R/{({I_{n - 1}},(\underline{x}'))}$
or the cokernel of the map
${({I_{n - 1}\cap(\underline{x}')})}/{(\underline{x}'){I_{n-2}}}\overset{x_l}\longrightarrow{({I_n}\cap (\underline{x}'))}/{{(\underline{x}'){I_{n-1}}}},$
then\\
$C_{n}={({I_n} \cap (\underline{x}'))}/{((\underline{x}'){I_{n-1}}+ x_{l}(I_{n-1}\cap (\underline{x}')))}$.
Thus if
$\grade(\underline{x}^{*})\geq l-1$,
 then $C_{n}=0$ for all $n$ and
${(({I_n},(\underline{x}')):{x_l})}/{({I_{n - 1}},(\underline{x}'))}\cong{({I_{n }} \cap (\underline{x}))}/{(\underline{x}){I_{n-1}}}$.
\end{remark}

From now on, we will assume that the filtration  $\mathcal F=\lbrace I^{n} \rbrace _{n =0}^{\infty}$ is $I$-adic filtration. Let $J = ({x_1},...,{x_d})$, where ${x_1},...,{x_d}$ is a superficial sequence in $I$, i.e., $J$ is a  minimal reduction of $I$. For $i\leq d-1$, set $J_{i}=({x_1},...,{x_{i}})$ (with the convention $J_{i}=(0)$ if $i\leq 0$ ),  and we denote $h_i( \underline x ,\mathcal{F})$, $k_i( \underline x ,\mathcal{F})$ and $h'_i(\underline x',\mathcal{F})$  for $I$-adic filtration by $h_i$ ,$k_i$ and $h'_i$, respectively.
\\

\begin{proposition}\label{3.6}
Let $d\ge 2$ and $J$ be a minimal reduction of $I$ such that
 $\sum\limits_{i=2}^d {{{( - 1)}^{i - 2}}{h_i}=1}$. Then $\depth G(I)\ge d-2.$
\end{proposition}
\begin{proof}
 If $\sum\limits_{i = 2}^d {{{( - 1)}^{i - 2}}{h_i} = 1}$, then by Remark 3.3,
 $\sum\limits_{n = 0}^\infty\lambda{({{I^{n + 1}}}/{J{I^n}})}- {e_1}(I) = 1$ and so by \cite[ Theorem 3.1]{W} we have $\depth G(I)\ge d-2.$
 \end{proof}

\begin{proposition}\label{3.7}
 Let $J$ be a minimal reduction of $I$ such that\\
 $\sum\limits_{n=1}^\infty\lambda({({{I^{n+1}} \cap J_{d-1}})/({J{I^n}\cap J_{d-1}}))=0}$. Then $\depth G(I)\ge{d-1}.$
\end{proposition}

\begin{proof}
By using induction on $n$, we prove that ${I^{n + 1}} \cap J_{d-1} = J_{d-1}{I^n}$ for every $n\ge0$. For $n =0$, there is nothing to prove. Assume that $n\ge1$ and ${I^{n }} \cap J_{d-1}=J_{d-1}{I^{n-1}}$. From the following equalities:
\[\begin{array}{ll}
{I^{n + 1}} \cap J_{d-1} &= J{I^n} \cap J_{d-1} \\&= (J_{d-1}{I^n} + {x_d}{I^n}) \cap J_{d-1}\\& = J_{d-1}{I^n} + ({x_d}{I^n} \cap J_{d-1})\\&= J_{d-1}{I^n} + {x_d}({I^n}\cap J_{d-1})\\& = J_{d-1}{I^n} + {x_d}(J_{d-1}{I^{n - 1}})=J_{d-1}{I^n}.
\end{array}\]
The fourth equality follows from the fact that $x_1,x_2, . . , x_d$ form a  regular sequence. Then by
 using Valabrega and Valla's theorem (see also \cite[Theorem 5.16]{V}) we have $x_{1}^{*}, x_{2}^{*}, . . . ,x_{d-1}^{*}$ are regular sequence and $\depth(G(I)) \geq d-1$.

\end{proof}

\begin{proposition}
Let $J$ be a minimal reduction of $I$ such that\\
 $\sum\limits_{n =1}^\infty \lambda(({{I^{n + 1}}\cap J_{d-1}})/({J{I^n} \cap J_{d-1}}))=1. $
 Then
 $\depth G(I)\ge {d-2}.$
\end{proposition}

\begin{proof}
Let $B_{n+1}$ be the kernel of the map
${H_1}({C_.}(n+1)) \to {H_{0}}({C_.}^\prime (n ))$
given in $(*)$. Consider the following exact sequence
\[0 \longrightarrow {H_d}({C_.}(n+1)) \longrightarrow ... \longrightarrow {H_1}({C_.}^\prime (n))\overset{x_d}\longrightarrow {H_1}({C_.}^\prime (n+1)) \longrightarrow {B_{n+1}} \longrightarrow 0.\]
Therefore by \cite[Lemma 3.2]{Hm} we have the following exact sequence
$$...\longrightarrow I^{n}\cap J_{d-1}/{J_{d-1}{I^{n-1}}}\overset{x_d}\longrightarrow {I^{n+1}\cap J_{d-1}/{J_{d-1}{I^{n }}}} \longrightarrow {B_{n+1}} \longrightarrow 0$$
and
$B_{n+1}={{I^{n + 1}} \cap J_{d-1}}/{J{I^n}\cap J_{d-1}}.$
Thus
\[\sum\limits_{n = 1}^\infty\lambda({{B_{n+1}}})=\sum\limits_{n = 1}^\infty\lambda(({{{I^{n + 1}} \cap J_{d-1}})/({{J{I^n}\cap J_{d-1}}}))=\sum\limits_{i=1}^d {{{( -1)}^{i-2}}}}{h_i}=1\]
and so by Proposition 3.6 we have
$\depth G(I)\ge d-2$.
\end{proof}

\begin{remark}\label{3.9}
 Let $B_n$ be the kernel of the map ${H_1}({C_.}(n))\to {H_{0}}({C_.}^\prime (n -1))$ given in $(*)$. Then for each $n$ we have the exact sequence
\[0\to {H_d}({C_.}(n)) \to ... \to {H_1}({C_.}^\prime (n - 1))\mathop  \to \limits^{{x_d}} {H_1}({C_.}^\prime (n)) \to {B_n} \to 0.\]

Therefore, for each $n$ we have
\[\lambda ({B_n}) = \sum\limits_{i=2}^d {{{( - 1)}^{ i - 2}}\lambda ({H_i}(} {C_.}(n))) + \sum\limits_{i = 1}^d {{{( - 1)}^{i-1}}\Delta [\lambda ({H_i}({C_.}^\prime } (n)))].\]
 By using \cite[Lemma 2.7]{H1}
\[\sum\limits_{n = 2}^\infty  {(n-1)\Delta [\lambda ({H_i}({{C'}_.}(n)))]}=-\sum\limits_{n=1}^\infty \lambda ({H_i}({{C'}_.}(n)))=- h'_i\]
and
\[\sum\limits_{n = 1}^\infty  {\Delta [\lambda ({H_i}({{C'}_.}(n))]}= 0.\]
Since
$\sum\limits_{n = 2}^\infty  {(n - 1)\lambda ({B_n}) = \sum\limits_{i = 2}^d {{{( - 1)}^{i - 2}}{k_i}} }-\sum\limits_{i = 1}^d {{{( - 1)}^{i - 1}}{{h'}_i}}$
and
$\sum\limits_{n = 1}^\infty {\lambda ({B_n})}$ is equal to ${\sum\limits_{i = 2}^d {{{( - 1)}^{i - 2}}{h_i}}}$,
 then by Lemma \ref{3.1},
 $\sum\limits_{i = 1}^d {{{( - 1)}^{i - 1}}{{h'}_i}}=0$ if and only if $\depth G(I)\geq d-1 $ if and only if $\sum\limits_{n=2}^\infty {(n-1)\lambda ({B_n})}=0.$
Thus
$\sum\limits_{i = 2}^d {{{( - 1)}^{i - 2}}{k_i}} = 1$ cannot be happen.
Now assume that $I$ is integrally closed ideal.  Since $I^{2}\cap J=JI$ (see \cite[Lemma 2.1]{Mn}), then by \cite[Lemma 3.2]{Hm} we have $H_{1}(C.(2))=0$ and $B_{2}=0$.
So if $I$ is integrally closed ideal, then $\sum\limits_{i = 2}^d {{{( - 1)}^{i - 2}}{k_i}}=2$ cannot  happen, because in this case
$\lambda(B_{2})=1$
and this is a contradiction.
 So we have the following three cases:
\item[(1)]
If $\sum\limits_{i = 2}^d {{{( - 1)}^{i - 2}}{k_i}}= 2 $, then we have $\sum\limits_{i = 1}^d {{{( - 1)}^{i - 1}}{{h'}_i}}=1$ and  $\sum\limits_{n = 2}^\infty  {(n - 1)\lambda ({B_n})}=1$. In this case we obtain $\lambda(B_{2})=1$ and $\lambda(B_{n})=0$ for any $n \ne 2$. Therefore $\sum\limits_{n = 1}^\infty  {\lambda ({B_n}) = \sum\limits_{i = 2}^d {{{( - 1)}^{i - 2}}{h_i}} }=1 $ and by Proposition \ref{3.6} $\depth G(I)\geq d-2$.
\item[(2)]
If we assume $I$ to be integrally closed ideal and $\sum\limits_{i =2}^d {{{( - 1)}^{i-2}}{k_i}}=3, $
then we have
$\sum\limits_{n = 2}^\infty (n - 1)\lambda ({B_n})=1 ~or ~2$. The case
$\sum\limits_{n = 2}^\infty (n - 1)\lambda ({B_n})=1$
cannot happen because
$\lambda(B_{2})=0$.
If  $\sum\limits_{n = 2}^\infty (n-1)\lambda ({B_n})=2$, then $\lambda(B_{3})=1$ and $\lambda(B_{n})=0$ for all $n\ne 3$,
so
$\sum\limits_{i = 2}^d {{{( - 1)}^{i - 2}}{h_i}}=1$
and
$\depth G(I)\geq d-2$.
\item[(3)]
If we assume $I$ to be integrally closed ideal and $\sum\limits_{i = 2}^d {{{( - 1)}^{i - 2}}{k_i}}= 4$, then we have $\sum\limits_{n = 2}^\infty(n - 1)\lambda ({B_n})=2~or~3$.
If $\sum\limits_{n = 2}^\infty (n - 1)\lambda ({B_n})=2$, then $\depth G(I)\geq d-2$.
  If $\sum\limits_{n = 2}^\infty(n - 1)\lambda ({B_n})=3$, then $\lambda(B_{4})=1$ and $\lambda(B_{n})=0$ for any $n\ne 4$. Thus $\sum\limits_{i = 2}^d {{{( - 1)}^{i - 2}}{h_i}}=1$ and $\depth G(I)\geq d-2$ .

\end{remark}

\noindent
In the following result we compare \cite[Theorem 3.1 and 3.3]{Cpr} with \cite[Theorem 3.1]{W}.
\begin{theorem}
Let $J$ be a minimal reduction of $I$. Suppose that one of the following
conditions holds:
\begin{enumerate}
\item ${e_2}(I)=\sum\limits_{n = 2}^\infty  {(n - 1)\lambda({{I^{n }}}/{J{I^{n-1}}})}-2; $
\item $I$ is integrally closed and ${e_2}(I)=\sum\limits_{n = 2}^\infty  {(n - 1)\lambda({{I^{n }}}/{J{I^{n-1}}})}-i$, where $i=3,4$.
\end{enumerate}
Then  ${e_1}(I)=\sum\limits_{n = 1}^\infty  {\lambda({{I^{n }}}/{J{I^{n-1}}})}-1.$
\end{theorem}
\begin{proof}
$(1)$ If ${e_2}(I)= \sum\limits_{n = 2}^\infty{(n - 1)\lambda({{I^{n }}}/{J{I^{n-1}}})}-2 $, then $\sum\limits_{i = 2}^d {{{( - 1)}^{i - 2}}{k_i}} = 2   $ and by Remark 3.9 $\sum\limits_{n = 1}^\infty  {\lambda ({B_n}) = \sum\limits_{i = 2}^d {{{( - 1)}^{i - 2}}{h_i}} }=1 $.  Therefore by Remark 3.3 we have ${e_1}(I)=\sum\limits_{n = 1}^\infty  \lambda{({{I^{n }}}/{J{I^{n-1}}})}-1$.\\

$(2)$ If I is  integrally closed and ${e_2}(I)= \sum\limits_{n = 2}^\infty{(n - 1)\lambda({{I^{n }}}/{{J{I^{n-1}}}})}-i$ where $i=3,4$ then $\sum\limits_{i = 2}^d {{{( - 1)}^{i - 2}}{k_i}} = 3~~or~4$. Therefore by Remark 3.9 $\sum\limits_{n = 1}^\infty{\lambda ({B_n})=\sum\limits_{i = 2}^d {{{( - 1)}^{i - 2}}{h_i}} }=1 $  and so ${e_1}(I)=\sum\limits_{n = 1}^\infty \lambda{({{I^{n}}}/{{J{I^{n-1}}}})}-1.$

\end{proof}
The following example shows that the converse of Theorem 3.10 in general is not true.
\begin{example}
Let $R=k[x,y]_{(x,y)}$, where $k$ is a field and $I=(x^6,y^6,x^5y+x^2y^4)$. Then by using Macaulay 2 \cite{Gs} we  obtain the following Hilbert polynomial
\[{P_I}(n) = 36\left( {\begin{array}{*{20}{c}}
{n + 1}\\
2
\end{array}} \right) - 15\left( {\begin{array}{*{20}{c}}
n\\
1
\end{array}} \right) + 11\]
and ${e_1}(I)= \sum\limits_{n = 1}^\infty  \lambda{({{I^{n }}}}/{{J{I^{n-1}}})}-1$ but ${e_2}(I)=\sum\limits_{n = 2}^\infty {(n - 1)\lambda({{I^{n }}}/{J{I^{n-1}}})}-3 $ and $I$ is not integrally closed.
\end{example}

\section{The depth of associated graded ring and the reduction number}
Let $A=\bigoplus_{m\geq 0}A_m$ be a Notherian graded ring where $A_0$ is an Artinian local ring, $A$ is generated by $A_1$ over $A_0$ and $A_{+}=\bigoplus_{m>0}A_m$. Let $H_{A_{+}}^i(A)$ denote the i-th local cohomology module of $A$ with respect to the graded ideal $A_+$ and set $a_i(A)=\max\{m\vert\ \  [H_{A_{+}}^i(A)]_m\neq 0\}$ with the convention $a_i(A)=-\infty$, if $H_{A_{+}}^i(A)=0$. The Castelnuovo-Mumford regularity is defined by $\reg (A):=\max\{a_i(A)+i\vert\ \ i\geq 0\}.$ In the following theorem, for simplicity, we use $a_i$ instead of $a_i(G(I))$.

\begin{theorem}
Let $J$ denote a minimal reduction of $I$. Suppose that one of the following
conditions holds:
\begin{enumerate}
\item ${e_2}(I)  = \sum\limits_{n = 2}^\infty{(n - 1)\lambda({{I^{n }}}/{J{I^{n-1}}})}-2;$
\item I is integrally closed and ${e_2}(I)= \sum\limits_{n=2}^\infty{(n - 1)\lambda({{I^{n }}}/{J{I^{n-1}}})}-3$.
\end{enumerate}
Then  $r(I)$ is independent.
\end{theorem}
\begin{proof}
(1) By Remark 3.9,  $\sum\limits_{i = 2}^d {{{( - 1)}^{i - 2}}{k_i}}= 2$ and $\depth G(I)\ge d-2$. If $\depth G(I)\ge d-1$, then by \cite[Theorem 2]{M1} $r_{J}(I)= \reg(G(I))$ and so the result in this case follows. Now we assume that $\depth G(I)= d-2$ and $\sum\limits_{i = 2}^d {{{( - 1)}^{i - 2}}{k_i}} = 2$. By Remark 3.9, $\lambda(B_{2})=1$ and $\lambda(B_{n})=0$ for any $n\ne 2$. Therefore we have $\lambda ({{I^2} \cap {J_{d - 1}}}/{{J_{d - 1}}I})= 1$, $\lambda ({{I^n} \cap {J_{d - 1}}}/{{{J_{d - 1}}{I^{n - 1}}}}) = 0$ for any $n\ne 2$, \\
$\lambda ({({I^2} + {J_{d-2}}):{x_{d - 1}}}/{I}) = 1$ and $\lambda ({({I^n} + {J_{d - 2}}):{x_{d - 1}}}/{{{I^{n - 1}} + {J_{d - 2}}}})= 0$ for any $n\ne 2$. If $\depth G(I)=d-2$, then by applying \cite[Theorem 2.1]{M} there are two cases:\\
(i) If $a_{d-1}\leq a_{d}$, then $r_{J}(I)=a_{d}+d=reg(G(I))$.\\
(ii)If $a_{d}<a_{d-1}$, then $r_{J}(I)\leq a_{d-1}+d-1=reg(G(I))$ and ${a_{d - 1}} = \max \{ n|{(I^{n + d-1}} + {J_{d - 2})}:{x_{d - 1}} \ne {I^{n + d - 2}} + {J_{d - 2}}\} $. Therefore by the above process ${a_{d - 1}}=3-d$ and so for any reduction $J$ of $I$, $r_{J}(I) \leq 2$. Hence $r(I)$ is independent.\\
 (2) Let $i=3$. Then by Remark 3.9, $\lambda(B_{3})=1$ and $\lambda(B_{n})=0$ for any $n\ne 3$. Thus $$\lambda ({{I^3} \cap {J_{d - 1}}}/{{{J_{d - 1}}I^2}}) = 1
 =\lambda ({({I^3} + {J_{d - 2}}):{x_{d - 1}}}/{I^2})$$ and $$\lambda ({{I^n} \cap {J_{d - 1}}}/{{{J_{d - 1}}{I^{n - 1}}}}) = 0=\lambda ({({I^n} + {J_{d - 2})}:{x_{d - 1}}}/{{{I^{n - 1}} + {J_{d - 2}}}})$$ for any $n\ne 3$. If $\depth G(I) = d-2$, then by \cite [Theorem 2.1]{M} there are two cases:\\
(i) If $a_{d-1}\leq a_{d}$, then $r_{J}(I)=a_{d}+d=reg(G(I))$.\\
(ii)
If $a_{d}<a_{d-1}$, then $r_{J}(I)\leq a_{d-1}+d-1=reg(G(I))$ and ${a_{d - 1}}=\max \{ n|{(I^{n + d-1}} + {J_{d - 2})}:{x_{d - 1}} \ne {I^{n + d - 2}} + {J_{d - 2}}\} $. Therefore by the above process ${a_{d - 1}}=4-d$ and so for any reduction $J$ of $I$, $r_{J}(I) \leq 3$. Since $I$ is integrally closed, then $r(I)$ is independent.
\end{proof}

Let $\fa$ be an ideal of grade at least $1$ in a Noetherian ring $R$. The Ratliff-Rush closure of $\fa$ is defined as the ideal $\widetilde{\fa}=\cup_{n\geq 1}({\fa}^{n+1}:{\fa}^n).$ It is a refinement of the integral closure of $\fa$ and $\widetilde{\fa}=\fa$ if $\fa$ is integrally closed (see \cite{Rr}).

\begin{proposition}(compare with \cite[Theorem 3.3(b)]{Cpr})
  Let $d=3$ and $J$ be a minimal reduction of $I$. If $\widetilde{I}=I$ and
 ${e_2}(I)=\sum\limits_{n = 2}^\infty  {(n - 1) \lambda({{I^{n }}}/{J{I^{n-1}}})}-3 $
then  $\depth G(I)\ge 1.$
\end{proposition}

\begin{proof}
If
${e_2}(I)=\sum\limits_{n = 2}^\infty  {(n -1)\lambda({{I^{n }}}/{J{I^{n-1}}})}-3,$
then
$\sum\limits_{i = 2}^d {{{( - 1)}^{i - 2}}{k_i}} = 3$.
Hence by Remark 3.9,
$\sum\limits_{i = 1}^d {{{( - 1)}^{i - 1}}{{h'}_i}}=1$
and
$\sum\limits_{n = 2}^\infty (n - 1)\lambda ({B_n})=2.$
Therefore we have the following cases:
\begin{enumerate}
\item[1)] $\lambda (B_{3})=1$ and $\lambda (B_{n})=0$  for any $n \ne 3$ .\\
or
\item[2)]$\lambda (B_{2})=2$ and  $\lambda (B_{n})=0$ for any $n \ne 2$.

\end{enumerate}
For the first case, by Proposition 3.6.  $\depth G(I)\ge 1$ and the result follows. Let us consider the second case. Since $\widetilde{I}=I$, we obtain that $I^2 : x = I$ for every superficial element $x$. Since $\lambda (B_{2})=2$, it follows $\lambda (I^2 \cap J_{2}/J_{2}I)=2$. It therefore follows $\lambda (I^2 \cap J_{1}/J_{1}I)=0$. Hence
$\sum\limits_{i = 1}^d {{{( - 1)}^{i - 1}}{{h'}_i}} \ge 2$ and this is a contradiction.

\end{proof}

Northcott in \cite{N} proved that
${e_1(I)} \ge {e_0(I)} - \lambda({R}/{I})$ and after that Huneke in \cite{Hu} showed that ${e_1(I)} = {e_0(I)}-\lambda ({R}/{I})$ if and only if ${I^2}= JI$. When this is the case, $G(I)$ is Cohen-Macaulay.
 Sally in \cite{S} proved that if $d\geq 2$, ${e_1(I)}-{e_0(I)} + \lambda ({R}/{I}) = 1$ and $e_{2}(I)\ne 0$, then  $\depth G(I)\ge d-1$ (see also
 \cite[Corollary 4.15]{Hm} and \cite[Proposition 3.1]{Gr1}). Also Itoh in \cite{I} with this conditions proved that if $I$ is integrally closed, then $G(I)$ is Cohen-Macaulay.
In the following example we show that the integrally closedness of $I$ is essential for the Cohen-Macaulayness of $G(I)$.

The following example appears in \cite{Cpr}.
\begin{example}
Let $R = k[x,y,z]_{(x,y,z)}$, where $k$ is a field, and $I=(x^2-y^2,y^2-z^2,xy,yz,xz)$. Then, by Macaulay 2, we have
${e_0}(I)=8$, ${e_1}(I)=4$, ${e_2}(I)=0$ and
$\depth G(I) =0.$
\end{example}

\begin{lemma}
Let $J$ be a minimal reduction of $I$. If ${e_1(I)}-{e_0(I)}+\lambda ({R}/{I}) = 2$ and $I=\bar I$,
then  $\depth G(I)\ge d-1. $
\end{lemma}

\begin{proof}
By using the Sally machine and the good behaviour of $e_1(I)$ modulo superficial elements (see \cite[Proposition 1.2]{Erv} and \cite[Lemma 2.2]{Hm}), we may reduce the statement to dimension two. If ${e_1(I)} - {e_0(I)} + \lambda ({R}/{I})=2$, by \cite[Corollary 1.5]{R} we have  ${r_J}(I) \le 3$ for any minimal reduction $J$ of $I$. If there exist a minimal reduction $J$ of $I$ such that ${r_J}(I) = 2$,  then by \cite[Lemma 2.1]{Mn} $G(I)$ is Cohen-Macaulay.
If for any minimal reduction $J$ of $I$, ${r_J}(I) = 3$ then by \cite[Remark 1.8]{R} ${r_J}(I) \le {e_1(I)} - {e_0(I)} + \lambda ({R}/{I})+2-\lambda ({{I^2}}/{{JI}})$. Hence $r(I)=3\le 4-\lambda ({{I^2}}/{JI})$ and $\lambda ({{I^2}}/{JI}) \le 1$ and so by \cite[Corollary 1.7]{R}
$\depth G(I)\ge 1 $.
\end{proof}

\begin{proposition}
Let $J$ be a minimal reduction of $I$. If ${e_1(I)}-{e_0(I)}+\lambda ({R}/{I})\leq 3$ and $I=\bar I$,
then $r(I)$ is independent.
\end{proposition}

\begin{proof}
By Lemma 3.15 and the above explanation, we can assume ${e_1(I)}-{e_0(I)}+\lambda ({R}/{I})=3$ and also by \cite[Corollary 4.7]{Or} and \cite[Lemma 1.1]{H} we may assume that $d=2$. If $\depth G(I)\ge 1 $, then by \cite[Corollary 2.2]{M} we have $r_J(I)=\reg(G(I))$ and so $r(I)$ is independent. Now we may assume that $\depth G(I)=0$.
Since ${e_1(I)}-{e_0(I)}+\lambda ({R}/{I})=3$, by \cite[Corollary 1.5]{R} ${r_J}(I)\le 4$ for any minimal reduction $J$ of $I$. If for some $J$,  ${r_J}(I)=2$ then by \cite[Lemma 2.1]{Mn} $G(I)$ is Cohen-Macaulay and this is a contradiction with $\depth G(I)=0$.
 If for some $J$,  ${r_J}(I)= 4$ then by \cite[Remark 1.8]{R} ${r_J}(I)\le {e_1(I)}-{e_0(I)}+\lambda({R}/{I})+2 -\lambda ({{I^2}}/{{JI}})$. Therefore $\lambda ({{I^2}}/{{JI}})\le 1$ and so
$\depth G(I)\ge 1$ and this is also a contradiction.
Hence we can assume that for any minimal reduction $J$ of $I$, ${r_J}(I)=3$ and so $r(I)$ is independent.
\end{proof}

\begin{corollary}
 Let $J$ be a minimal reduction of $I$. If ${e_1(I)}-{e_0(I)}+\lambda ({R}/{I})\leq{r_J}(I)-1$ , $I =\bar I$ and $\depth G(I)\ge d-2$,
then  $r(I)$ is independent.

\end{corollary}

\begin{proof}
By using the Sally machine and the good behaviour of $e_1(I)$ modulo superficial elements, we may reduce the statement to dimension two. If ${e_1(I)} - {e_0(I)}+\lambda ({R}/{I})\leq{r_J}(I)-1 $, then by \cite[Remark 1.8]{R} ${r_J}(I)\le {e_1(I)}-{e_0(I)}+\lambda ({R}/{I})+2-\lambda ({{I^2}}/{{JI}})$. Hence $\lambda ({{I^2}}/{{JI}}) \le 1$ and so $\depth G(I)\ge d-1$. Therefore the result follows by \cite[Theorem 2]{M1}.
\end{proof}

\begin{example}
Let $R = k[x,y]_{(x ,y)}$ where $k$ is a field, and $I=(x^6,y^6,x^5y,x^3y^3,x^2y^4,xy^5)$. Then by Macaulay 2 we have
${e_0}(I)  =36$ , ${e_1}(I)  =15$ and $\lambda(R/I)=22$. Hence ${e_1}(I)-{e_0}(I)+\lambda(R/I)=1$ but
 $\depth G(I)=0$ and $r_{J}(I)=2=reg(G(I))$ for all minimal reduction $J$ of $I$.
\end{example}
\begin{example}
Let $R = k[x,y]_{(x ,y)}$ where $k$ is a field, and $I=(x^6,y^6,x^5y,x^3y^3,x^2y^4)$. Then by Macaulay 2 we have
${e_0}(I)  =36$ , ${e_1}(I)=15$ and $\lambda(R/I)=23$. Hence ${e_1}(I)-{e_0}(I)+\lambda(R/I)=2$ but
 $\depth G(I)=0$ and $r_{J}(I)=2=reg(G(I))$ for all minimal reduction $J$ of $I$.
\end{example}
\begin{example}
Let $R = k[x,y]_{(x ,y)}$ where $k$ is a field, and $I=(x^6,y^6,x^5y,x^3y^3,xy^5)$. Then by Macaulay 2 we have
${e_0}(I)=36$, ${e_1}(I)=15$ and $\lambda(R/I)=23$. Hence ${e_1}(I)-{e_0}(I)+\lambda(R/I)=3$ but
 $\depth G(I)=0$ and $r_{J}(I)=2=reg(G(I))$ for all minimal reduction $J$ of $I$.
\end{example}

\begin{example}
Let $R = k[x,y]_{(x ,y)}$ where $k$ is a field, and $I=(x^6,y^6,x^5y,x^2y^4,xy^5)$. Then by Macaulay 2 we have
${e_0}(I)  =36$, ${e_1}(I)=15$ and $\lambda(R/I)=24$. Hence ${e_1}(I)-{e_0}(I)+\lambda(R/I)=3$ but for two minimal reduction $J_1=(x^6,x^5y+y^6)$ and $J_2=(x^6,y^6)$ we have $r_{J_1}(I)=2$ and $r_{J_2}(I)=3$ and
$\depth G(I) =0 $ because $I$ is not integrally closed.
\end{example}
The following example is due to Huckaba and Huneke \cite{Hh}.
\begin{example}
Let $R = k[x,y,z]_{(x ,y,z)}$ where $k$ is a field of characteristic $\ne 3$. Let $\fa=(x^4,x(y^3+z^3),y(y^3+z^3),z(y^3+z^3))$ and set $I=\fa+{\fm}^5$ where $\fm$ is the maximal ideal of $R$. The ideal $I$ is a integral closer $m$-primary ideal whose associated
graded ring $gr_I (R)$ has depth 2. We checked that
${e_0}(I)  =76$ , ${e_1}(I)  =48$ and $\lambda(R/I)=31$  so ${e_1}(I)-{e_0}(I)+\lambda(R/I)=3$ and $r(I)$ is independent.

\end{example}
\begin{acknowledgement} 
We would like to thank deeply grateful to the referee for the careful reading
of the manuscript and the helpful suggestions
\end{acknowledgement}

%%%%%%%%%%%%%%%%%%%%%%%%%%%%%%%%%%%%%%%%%%%%%%%%%%%%%%%%%%%%%%%%%%%%%%%%%%%%%


\begin{thebibliography}{99}
\bibitem{Bh}
W. Bruns and J. Herzog, {\it Cohen-Macaulay rings}, Cambridge University Press, Cambridge,
UK, (1998).

\bibitem{Cpr}
A. Corso , C. Polini and M. E. Rossi , {\it Depth of associated graded rings via Hilbert coeffcients of
ideals}, J. Pure and Appl. Algebra, {\bf 201} (2005), 126-141.
\bibitem{Erv}
J. Elias, M. E. Rossi and G. Valla, {\it On the coefficients of the Hilbert polynomial}, J. Pure and Appl. Algebra, {\bf 108} (1996), 35-60
\bibitem{Gs}
D. R. Grayson and M. E. Stillman, {\it Macaulay 2, a software system for research in algebraic
geometry}, {Available at http:// www.math.uiuc.edu/Macaulay2}.

\bibitem{G1}
A. Guerrieri, {\it On the depth of the associated graded ring of an m-primary ideal of a Cohen-Macaulay local ring}, J. Algebra, {\bf 167} (1994), 745-757.
\bibitem{G2}
A. Guerrieri, {\it On the depth of the associated graded ring}, Proc. Amer. Math. Soc., {\bf 123} (1995), 11-20.
\bibitem{Gr1}
A. Guerrieri and M. E. Rossi, {\it Hilbert coefficients of Hilbert filtrations}, J. Algebra, {\bf 199} (1998), 40-61.
\bibitem{Gr}
A. Guerrieri and M. E. Rossi, {\it Estimates on the depth of the associated graded ring}, J. Algebra, {\bf 211} (1999), 457-471.
\bibitem{Ho}
L. T. Hoa, {\it Reduction numbers and rees algebras of powers of an ideal}, Proc. Amer. Math. Soc., {\bf 119} (1993), 415-422.
\bibitem{Ho1}
L. T. Hoa, {\it Reduction numbers of equimultiple ideals}, J. Pure and Appl. Algebra, {\bf 109} (1996), 111-126.

\bibitem{H}
S. Huckaba, {\it Reduction numbers of ideals of higher analytic spread}, Proc. Camb. Phil. Soc.,
{\bf 102} (1987), 49-57.

\bibitem{H1}
 S. Huckaba, {\it A d-dimensional extension of a lemma of Huneke's and formulas for the Hilbert coefficients}, Proc. Amer. Math. Soc., {\bf 124} (1996), 1393-1401.
\bibitem{Hh}
S. Huckaba and C. Huneke, {\it Normal ideals in regular rings}, J. Reine Angew. Math., {\bf 510} (1999), 63-82.
\bibitem{Hm}
S. Huckaba and T. Marley, {\it Hilbert coefficients and the depths of associated graded rings}, J. London Math. Soc., {\bf 56} (1997), 64-76.
\bibitem{Hu}
C. Huneke, {\it Hilbert functions and symbolic powers}, Michigan Math. J.,{\bf 34} (1987), 293-318.
\bibitem{I}
S. Itoh, {\it Hilbert coefficients of integrally closed ideals}, J. Algebra, {\bf 176} (1995), 638-652.
\bibitem{Mn}
A. Mafi and D. Naderi, {\it A note on reduction numbers and Hilbert-Samuel functions of ideals over Cohen-Macaulay rings}, Turkish J. Math., {\bf 40} (2016), 766-769.
\bibitem{M1}
T. Marley, {\it The coefficients of the Hilbert polynomial and the reduction number of an ideal}, J. London Math. Soc., {\bf 40} (1989), 1-8.

\bibitem{M}
T. Marley, {\it The reduction number of an ideal and the local cohomology of the associated graded ring}, Proc. Amer. Math. Soc., {\bf 117} (1993), 335-341.

 \bibitem{N}
D. G. Northcott, A note on the coefficients of the abstract Hilbert function, J. London Math.
Soc., {\bf 35} (1960) 209-214.
\bibitem{Nr}
D. G. Northcott and D. Rees, Reduction of ideals in local rings, Proc. Cambridge Philos. Soc.,
{\bf 50} (1954), 145-158.
\bibitem{Or}
K. Ozeki and M. E. Rossi, {\it The structure of the Sally module of integrally closed ideals}, Nagoya Math. J., {\bf 227}  (2017), 49-76.

\bibitem{P}
C. Polini, {\it A filtration of the Sally module and the associated graded ring of an ideal}, Comm. Algebra, {\bf 28} (2000), 1335-1341.
\bibitem{Rr}
L. J. Ratliff and D. Rush, {\it Two notes on reductions of ideals}, Indiana Univ. Math. J., {\bf 27} (1978), 929-934.

\bibitem{R}
M. E. Rossi, {\it A bound on the reduction number of a primary ideal}, Proc. Amer. Math. Soc., {\bf 128} (1999), 1325-1332.

\bibitem{Rv}
M. E. Rossi and G. Valla, {\it Hilbert functions of filtered modules}, Lecture Notes of the Unione Matematica Italiana, Vol. {\bf 9}, Springer-Verlag, Berlin; UMI, Bologna, 2010.

\bibitem{S1}
J. D. Sally, {\it Reductions, local cohomology and Hilbert functions of local ring}. In Commutative Algebra: Durham 1981. London Math. Soc. Lecture Notes Series no. {\bf 72}(Cambridge University Press, 1982), pp. 231-241.
\bibitem{S}
J. D. Sally, {\it Hilbert coefficients and reduction number 2}, J. Algebraic Geom., {\bf 1} (1992), 325-333.
\bibitem{St}
B. Strunk {\it Castelnuovo-Mumford regularity, postulation numbers, and reduction numbers}, J. Algebra, {\bf 311}(2007), 538-550.
\bibitem{T}
N. V. Trung, {\it Reduction exponent and dgree bound for the defining equations of graded rings},
Proc. Amer. Math. Soc., {\bf 101} (1987), 229-236.

\bibitem{Vv}
P. Valabrega and G. Valla, {\it Form rings and regular sequence}, Nagoya Math. J., {\bf 72} (1978), 93-101.

\bibitem{Va}
M. Vaz Pinto, {\it Hilbert functions and Sally modules}, J. Algebra, {\bf 192} (1997), 504-523.
\bibitem{V}
J. K. Verma, {\it Hilbert coefficients and depth of the associated graded ring of an ideal}, arXiv:0801.4866.
\bibitem{W}
H. J. Wang, {\it Hilbert coefficients and the associated graded rings}, Proc. Amer. Math. Soc., {\bf 128} (2000)
963-973.
\bibitem{W1}
H. J. Wang, {\it On the associated graded rings of ideals of reduction 2}, J. Algebra, {\bf 236} (2001), 192-215.
\bibitem{Wu}
Y. Wu, {\it Reduction numbers and Hilbert polynomials of ideals in higher dimensional Cohen-Macaulay local rings}, Math. Proc. Camb. Phil. Soc., {\bf 111} (1992), 47-56.
\end{thebibliography}
\end{document}